\documentclass[12pt]{article}

\usepackage{amsmath, amsthm, amssymb}

\newtheorem{teor}{Theorem}[section]

\newtheorem{remark}{Remark}[section]
\newtheorem{prop}{Proposition}[section]

\newtheorem{guess}{Conjecture}

\newcommand{\real}{\mathbb R}
\newcommand{\complex}{\mathbb C}

\newcommand{\K}{K\"{a}hler}

\begin{document}

\title{The Weierstrass representation for pluriminimal
submanifolds}
\date{}
\maketitle

\vspace{-2cm}
\begin{center}
{\large
Claudio Arezzo
$\!\!$\footnote{
Partially
supported by Cofin 2000: Geometria delle variet\`a
reali e complesse.}
\ \,\,\,
Gian Pietro Pirola $\!\!$\footnote{
Partially
supported by Cofin 1999: Spazi di moduli e teoria delle
rappresentazioni;
Far 2000 (Pavia): Variet\`{a} algebriche, calcolo
algebrico,
grafi orientati e topologici.}
\ \,\,\,
Margherita Solci $\!\!$\footnote{
Partially
supported by G.N.A.M.P.A.}
}
\end{center}

\pagestyle{myheadings}
\markboth{Pluriminimal submanifold}{Pluriminimal submanifold}
%\markright{}

{\it{1991 Math. Subject Classification:}} 58E12, 53A10.

\vspace{2cm}

{\small
{\bf Abstract}
\ In this note we prove a Weierstrass representation formula for pluriminimal 
submanifolds of euclidean spaces. 
We use this formula to produce new examples
of pluriminimal $4$-submanifolds of ${\real}^6$ and 
to prove that  
any affine algebraic manifold can be pluriminimally embedded 
into some euclidean
space in a non holomorphic manner.}

\section*{Introduction}

The classical methods to describe
minimal submanifolds of
riemannian manifolds
are complex analysis for $2-$dimensional
domains,
and the study of
the minimal equation
for hypersurfaces.
In the intermediate cases these
tools do not give a satisfactory
description of the picture. It is then natural to restrict
the class of minimal submanifolds.
For example, when
$M$ is a complex manifold of dimension $m$ (when
necessary we will denote by $J$ its complex structure),
$(X,g)$ is a riemannian manifold,
following Eschenburg and Tribuzy (\cite{et}) we set:

\medskip
\noindent {\bf Definition.}\
{\em An immersion $f\colon M \rightarrow X$
is called {\emph{pluriminimal}}
if the re\-stric\-tion to any smooth
complex curve in $M$
is a minimal immersion
into~$X.$}

\medskip
\noindent We remark that if $m=1$
pluriminimal is equivalent
to minimal.

\medskip

The first problem is to show that this class of submanifolds
contains interesting examples.

\medskip
In this paper we study the case when
$(X,g)$ is the euclidean space. We propose an analogue of the Weierstrass
representation for pluriminimal maps.
As for minimal surfaces, this formula
allows to construct
many examples, either by explicit
calculations, or by using
techniques of complex geometry to
establish existence results.
We give an application of both
ways by constructing infinite
families of pluriminimal immersions
in particular of ${\complex}^2$ into ${\real}^6$,
generalizing the example found by Furuhata
(\cite{fu}).

\medskip
We also solve
a general existence problem in
arbitrary dimension, once again
drawing analogy with the case of
minimal surfaces. In fact we prove
that all affine algebraic
manifolds (i.e. compact projective
minus an ample divisor) admit a
pluriminimal nonholomorphic immersion
into some euclidean space.
The corresponding result for minimal
surface in $3$-space has been proved
in \cite{pi}, where it is shown
that any compact Riemann
surface minus any set of
point can be minimally
immersed into ${\real}^3$.

\medskip
Using the Weierstrass representation
we also give a simple
proof of the fact that
pluriminimal immersions induce
a \K\ metric on the domain.
Thus these submanifolds can be seen as isometric
pluriharmonic immersions of
\K\ manifolds, which have been
extensively studied by many authors, in
particular we refer to the work of
Dajczer, Gromoll and Rodriguez
(\cite{dg2}, \cite{dg1},
\cite{dr2}, \cite{dr1}).
We underline that this relation
holds only for submanifolds
of euclidean spaces.
This suggests that pluriminimal immersions have a
variational characterization, which
greatly enhances interest in
their study, and which has been
succesfully used  by many authors to solve
rigidity questions for \K\ manifolds,
see e.g. Siu (\cite{si}), and
Jost-Zuo (\cite{jz}).

\medskip

Because of these considerations, it
seems natural to ask to which
extent the analogy with the two
dimensional case carries over.
In particular, we point out the
problem of the extension of
Osserman's Theorem (\cite{os}),
which states that, if the minimal
surface has finite total
curvature, the holomorphic
$1$-forms which appear in the
Weierstrass formula extend to meromorphic
data on a compact riemann surface.
We believe it would be very
interesting to find the geometric
hypothesis which allow to compactify
the pluriminimal submanifold
in such a way that the Weierstrass formula
extend to meromorphic data on the compactification.

\medskip

\noindent {\small
{\bf Acknowledgment.}
\ We wish to thank Professor Eschenburg for many clarifying conversations
on this topic.}
% Since all these considerations, it
% seems natural to ask to which extent
% the analogy with the two
% dimensional case carries over.
% In particular, we point out the
% problem of the extension for the
% Osserman Theorem (\cite{os}),
% finding the geometric hypothesis
% which allows to
% compactify the pluriminimal
% submanifold in such a way
% that the Weierstrass formula
% extend to meromorphic data
% on the compactification.

\section{The Weierstrass formula}
Let $f\colon M\rightarrow \mathbb R^n$ be a smooth map.
We can write:
\begin{equation}\nonumber
f(Q)=\int_{P}^{Q} df + f(P)=Re \int_{P}^{Q}(\omega_{1},
\dots,\omega_{n})+f(P),
\end{equation}
where the $\omega_{i}$ are $(1,0)$-forms on $M$,
that is they are smooth sections of the complex cotangent
bundle $\Omega^{1}_{M}$ of $M.$
In local coordinates:
$$\omega_{i}(z_{1},\dots,z_{m})=\sum_{j=1}^{m}\omega_{ij}
(z_{1},\dots,z_{m})dz_{j}.$$

\medskip
The {\em conformality tensor} will be the section of
$Sym^{2}\Omega^{1}_{M}$
defined by:
$$
\Omega = \sum _{i=1}^{n} \omega _i\otimes \omega _i.
$$

\medskip
\noindent
The following result characterizes
the pluriminimal immersions.

\begin{teor}
\label{weie}
Let $\omega_{1},\dots,\omega_{n}$ be $(1,0)$ smooth forms
of $M,$
such that
$Re\  \omega_{i}$ is exact for every $i.$
Then
\begin{equation}\label{f}
f(Q)=Re \int_{P}^{Q}(\omega_{1},
\dots,\omega_{n})+ const,
\end{equation}
defines a
pluriminimal immersion
if and only if:
\begin{enumerate}
\item[a)] the $\omega_{i}$ are closed holomorphic;
\item[b)] \label{co}
the conformality tensor
vanishes:
\begin{equation}\label{coeq}
\sum _{i=1}^{n} \omega _i\otimes \omega _i=0;
\end{equation}
\item[c)]
\label{nozero}
the (complex) jacobian matrix
$(\omega_{ik})$ has maximal rank at every
point.
\end{enumerate}
\end{teor}

% , or equivalently, having set $W = span_{{\complex}}
% \{\omega_{1},\dots,\omega_{n}\}, $ then the linear system
% $|\bigwedge^{m}W|$ is base point free;

\begin{proof}
The classical Weierstrass representation
formula for minimal surfaces
implies that
if the properties $a),$ $b)$ and $c)$ hold, the map
$f$ defined in (\ref{f}) is a pluriminimal immersion.

Conversely,
let us first prove that each
one of the $\omega_{i}$ is
holomorphic: indeed, we know that
$\omega_{i}\mid_{C}$ is holomorphic on
each holomorphic curve $C$
(see \cite{os}).
Chosen $P\in M,$ we can find
local coordinates $z_1,
\dots, z_m$ such that $z_{j}(P)=0$
for every $j.$
Since
$\omega_{i} = \sum _{j=1}^m \omega_{ij}(z_1,
\dots, z_m) dz_j$, we can write

$$\bar{\partial}\omega_{i}=\sum_{j,k=1}^{m}
\frac{\partial \omega_{ij}}{\partial
\bar{z}_k}d\bar{z}_k \wedge dz_j\,\, .$$

Restricting this form to each line $z=a\xi$, $a= (a_j) \in {\complex}^n$,
$\xi
\in {\complex}$, we get
$$
\sum_{j,k=1}^{m}\frac{\partial \omega_{ij}}
{\partial \bar{z}_k} a_j \bar{a}_k = 0,$$
which clearly implies $\frac{\partial \omega_{ij}}{\partial
\bar{z}_k} = 0$ for any $j$ and $k$.

Let us now prove that $\omega_{i}$ has
to be closed: we know that
$Re\ \omega_{i}$ is closed and
$\omega_{i}$ is holomorphic.
Then
$0= (\partial + \bar{\partial})(\omega _{i} +
\bar{\omega}_{i}) = \partial \omega_{i}
+  \bar{\partial}\bar{\omega}_{i}$.
Since $\partial  \omega_{i}$ is of type
$(2,0)$
and $\bar{\partial}\bar{\omega}_{i}$ is of type
$(0,2)$ we get $\partial
\omega_{i}=0$ and  $\bar{\partial}\bar{\omega}_{i}=0$ which
immediately imply
$d\omega_{i} = 0$.

The conformality condition $b)$ follows directly from the fact that
given any vector in complexified tangent space $v \in T_{{\complex}}M$,
there exists a  complex curve with $v$ as tangent vector. On this curve
$f$ has to be minimal, which, by the classical
Weierstrass representation
formula, implies $\Omega (v,v) = 0$.

Condition $c)$
follows
by contradiction.
Indeed, if
$v \in \ker D_{{\complex}}(f)$,
where $D$ stands for the jacobian,
we can take a complex curve $C$
in $M$ tangent to $v$. By
restricting $f$ to this Riemann surface we get $\omega_i (v) =0$ for any
$i$, and then $f\!\mid _{C}$ is not an immersion.
\end{proof}

The geometrical
meaning of condition
$c)$ is given in the following:

\medskip
%\noindent {\bf Remark.}
\begin{remark} \label{gauss}
{\rm Let $W = span_{{\complex}}
  \{\omega_{1},\dots,\omega_{n}\}$ be
the space generated by the
$\omega_{i},$
and consider the natural map $\lambda\colon
\bigwedge^{m}W\to H^{0}(M,\Omega^{m}_{M}),$ where
$\Omega^{m}_{M}$
is the canonical bundle of $M.$
Then,
the immersion
property
$c)$ holds if and only if
the linear system
$|\lambda(\bigwedge^{m}W)|$
is base point free.
We note that the associate map
$g\colon M\to |\lambda(\bigwedge^{m}W)|$
is the composition of the (complex) Gauss map
with the Pl\"ucker embedding.
This explains why it is more difficult
for $m>1$ to see the appearance
of the Gauss map in the Weierstrass formula.}
\end{remark}

\medskip
On the other hand, condition $b)$
can be used to give simple proofs of two results of
Dajczer-Rodriguez (\cite{dr2}).

\begin{prop}
The riemannian metric induced by a pluriminimal immersion is \K .
\end{prop}
\begin{proof}
Let us consider complex local coordinates ${z_k}$ on $M$, and consider
the real and imaginary parts as real coordinates $z_k = x_k + i y_{k}$.
In real notation we can write
$$\omega_ j = \sum _{k=1}^{m} \alpha _{jk}
dx_k - \beta_{jk}dy_{k} + i (\sum _{l=1}^{m}\alpha _{jl}
dy_{l} + \beta_{jl}dx_l)\,\, ,$$
$j=1, \dots , n$, and
in matrix form
\begin{equation}
\label{zu}
\underline{\omega} = A\underline{x} - B\underline{y} + i(A
\underline{y} + B\underline{x}).
\end{equation}

Clearly $df(\frac{\partial}{\partial x_j}) = (\alpha
_{1j},
\dots ,
\alpha _{nj})$ and $df(\frac{\partial}{\partial y_{j}}) = -(\beta
_{1j},
\dots ,
\beta _{nj})$ for $j=1, \dots , m$.
Therefore, having set $g=f^*(eucl) $,
$$g(\frac{\partial}{\partial x_i},\frac{\partial}{\partial x_j})=
\alpha _{1i}\alpha _{1j} + \dots + \alpha _{ni}\alpha _{nj}\,\, ,$$
$$g (\frac{\partial}{\partial x_i},\frac{\partial}{\partial
y_{j}})=
-\alpha _{1i}\beta _{1j} - \dots - \alpha _{ni}\beta _{nj}\,\, ,$$
$$g (\frac{\partial}{\partial y_{i}},\frac{\partial}{\partial
y_{j}})=
\beta _{1i}\beta _{1j} + \dots + \beta _{ni}\beta _{nj}\,\, .$$

Thus we can write the matrix associated to $g$ as
$$\begin{pmatrix}
A & -B
\end{pmatrix} \begin{pmatrix}
A^{t} \\
-B^{t}
\end{pmatrix}= \begin{pmatrix}
AA^{t} & -AB^{t} \\
-BA^{t} & BB^{t}
\end{pmatrix}.$$

On the other hand
the vanishing of the
tensor $\Omega = \sum _{i=1}^{n} \omega _i\otimes \omega _i = 0$
can be written as $\Omega = \underline{\omega} \otimes
\underline{\omega} =0$,
which gives, by equation \ref{zu}, the following system of equations

$$ -AB^{t} = -BA^{t} = 0 \, , \quad AA^{t} = BB^{t}\,\,\,.$$

Therefore the matrix associated to $g$ is hermitian
and its associated form can be
written, as in the classical case of minimal surfaces, as
$\sum _{r=1}^{n}\omega_r\wedge\bar{\omega}_r$ which is clearly positive
definite and of type $(1,1)$. It is also closed since each $\omega_r$ is
closed.
\end{proof}
\medskip

The above proposition is crucial to
link our definition to more standard notions in the theory
of higher dimensional submanifolds of euclidean spaces. In particular let
us observe that pluriminimal immersions are
part of a broader class of submanifolds studied in general by many authors
(e.g.
\cite{dg2}, \cite{dg1}, \cite{dr2}, \cite{dr1}).

\begin{prop}
The second fundamental form $B$ of a pluriminimal immersion satisfies
$$B(X, JY) = B(JX,Y)\,\, \, ,$$
i.e. a pluriminimal immersion is circular.
\end{prop}
\begin{proof}
Since the induced metric is \K , at every point of $M$ we can choose an
orthonormal basis for the tangent space of the form
$\{e_1, \dots , e_m, Je_1,\dots , Je_m\}$. Since the map
restricted to every holomorphic direction has to be minimal, $B(e_j, e_j)
+ B(Je_j, Je_j) = 0$ for any $j=1, \dots , m$. Therefore
$$B(e_j+e_k,e_j+e_k)+ B(J(e_j+e_k),J(e_j+e_k))=
2(B(e_j,e_k)+B(Je_j,Je_k))=0$$
Thus, $(B) = -(J^tBJ)$, which implies directly the conclusion.
\end{proof}

\bigskip

\section{Constructions of pluriminimal immersions}

We look for holomorphic functions of two
complex variables $x$ and $y$,
whose differentials satisfy the quadratic relation

\begin{equation}
\label{qqcc}
dP_{1}\otimes dP_{2}= dP_{3}\otimes dP_{4} + dP_{5}\otimes dP_{6}.
\end{equation}

By a diagonalization process
on the above tensor,
in such a way that the condition $b)$
of the theorem \ref{weie}
is satisfied, we can write the map $f$
defined in (\ref{f}) as:

\begin{equation}\label{im6}
(x,y)\mapsto
\begin{pmatrix}
Re (P_{1}+P_{2}) \\
Im (P_{1}-P_{2}) \\
Re (P_{3}+P_{4}) \\
Im (P_{4}-P_{3}) \\
Re (P_{5}+P_{6}) \\
Im (P_{6}-P_{5})
\end{pmatrix}
\end{equation}

Let $W=$ span$_{\complex}[dP_{1},\dots, dP_{6}];$
if $W$ satisfies the condition $c),$
we can choose, in local coordinates,
$P_{3}=x$ and $P_{4}=y.$ Moreover,
we set $P_{1}(x,y)=xy.$ Then the equation (\ref{qqcc})
translates into the system:
\begin{eqnarray}\label{esempio}
\left\{ \begin{array}{ll}
   y (P_{2})_{x}=(P_{4})_{x} \\
   x (P_{2})_{y}=(P_{6})_{y} \\
   y (P_{2})_{y}+x (P_{2})_{x}=
     (P_{4})_{y}+(P_{6})_{x}
\end{array}
\right.
\end{eqnarray}

A simple calculation shows that the solutions of (\ref{esempio})
are of the following form:

\begin{eqnarray}
&& P_{2}(x,y)=- \frac{g^{\prime}(x)+f^{\prime}(y)}{2}
    \nonumber \\
&& P_{4}(x,y)= f(y)-y \ \frac{g^{\prime}(x)+f^{\prime}(y)}{2}
    \nonumber \\
&& P_{6}(x,y)=g(x)-x \ \frac{g^{\prime}(x)+f^{\prime}(y)}{2}
    \nonumber
\end{eqnarray}
where $f$ and $g$ are arbitrary holomorphic functions of one complex
variable.

\smallskip
We observe that by taking $f,g$ entire functions
we obtain pluriminimal non holomorphic
immersions of $\mathbb C^{2}$ in $\mathbb R^{6}.$
The only example of this sort known to
us is due to Furuhata (\cite{fu}),
and belongs to our class for $f=x^3$ and $g=0,$
up to a real constant.
Nevertheless by direct computation it is
possible to show that Furuhata's
example is not an embedding, i.e.
the map is not injective. We believe
this should be true for all such maps.

\begin{guess}
Any complete pluriminimal immersion from
${\complex}^2$ to ${\real}^6$
is not an embedding.
\end{guess}

It is clear that with a similar procedure, choosing
meromorphic or algebraic functions,
we could construct families of pluriminimal immersions
of more comlicated domains.
For other examples, see also \cite{dg3}.
\medskip

We underline that the Weierstrass representation theorem
also allows to prove general existence theorems, without
explicitely finding
the holomorphic differentials,
in total analogy with the theory of
minimal surfaces in ${\real}^3$.

\medskip
We now prove that
every affine algebraic manifold $X$
admits a pluriminimal embedding into some euclidean
space.
Let us then start with a
smooth projective manifold, $M$, of complex dimension
$m$, and let $H$ be a hyperplane section.
Set $X=M\setminus H.$
By Hirzebruch-Riemann-Roch' s
Theorem (see, for example,
\cite{ha} and \cite{hi}) we know that,
for large $n:$
$$\dim H^0(M,{\mathcal{O}_{M}}(nH))=
\displaystyle\frac{H^m}{m!}n^{m}+P(n),$$
where $P$ is a polynomial of degree
at most $m-1$ in the variable $n.$
In order to construct
holomorphic $(1,0)$-forms
on the manifold $M\setminus H$,
we consider the image $V$ of the
exterior differential
$d \colon H^{0}(M,{\mathcal{O}_{M}}(nH))
\to H^{0}(M,\Omega^1_{M}((n+1)H)$.
Note that the forms in $H^{0}(M,\Omega^1_{M}
((n+1)H)$ are holomorphic on $X.$

By restriction, the cup product map
$$\mu_{n} \colon Sym^2H^0(M, \Omega^1_{M}((n+1)H))
\rightarrow H^0(M,Sym^2\Omega
^1_{M}(2(n+1)H))$$
defines an application
$$\mu^{\prime}_{n} \colon Sym^2 V \rightarrow
  H^0(M,Sym^2\Omega^1_{M}(2(n+1)H)). $$

Any element in the kernel of $\mu^{\prime}_{n}$ represents a quadratic
relation satisfied by holomorphic $(1,0)$-forms, and therefore
we can diagonalize the tensor
in order to satisfy the condition $b)$
in the theorem \ref{weie}.

We then estimate the dimension
of $\ker (\mu^{\prime}_{n}).$
Once again by Hirzebruch-Riemann-Roch,
for large $n$ we have:
$$\dim Sym^2H^0(M, \Omega^1_{M}((n+1)H))=
(\frac{m(m+1)}{2})\frac{H^m}{m!}(2n+2)^m+Q(n),$$
$Q$ being a polynomial of degree at most $m-1.$
Since $\dim Sym^2 V$ grows as $n^{2m},$
the map $\mu^{\prime}_{n}$
has nontrivial kernel for $n$
large enough.

\bigskip

At this point we can construct a
pluriminimal map by associating
to a nontrivial element $\gamma$ of
$\ker \mu^{\prime}_{n}$
a set of independent
exact $(1,0)-$forms
$dF_1, \dots , dF_k$,
where $k$ is the rank of $\gamma$,
satisfying $\sum _{j=1}^k dF_j \otimes dF_j = 0$.
Then, the map $\phi\colon M\setminus H \rightarrow {\real}^k$ defined by
$$\phi (p) = Re (F_1, \dots, F_k) + const$$
is a pluriminimal map.

As proved by Arezzo, Micallef and Pirola (\cite{amp}), the fact that the
kernel contains
nontrivial elements
easily implies that
$\phi$ is not holomorphic
w.r.t. any complex structure.

It is immediate to check that for $n$ big enough one can find
$F_{i}$ s.t. $\phi$ is an embedding.
Moreover, the previous estimates prove the
following:

\begin{teor}
Let $X$ an affine algebraic variety of dimension $m$;
then we can find an integer $k(m)$ such that
for every $n\geq k(m)$ there exist
pluriminimal embeddings
$$\phi \colon X \rightarrow {\real}^n$$
such that
\begin{enumerate}
\item the Gauss map $g$ defined in Remark
\ref{gauss} is algebraic;
\item the moduli number $c_{n}$
satisfies $$c_{n}=O(n^{2m})_{n\to\infty}.$$
\end{enumerate}
\end{teor}

\begin{remark}
{\rm We note that in the previous theorem
if $n$ is even $\phi$ is not holomorphic w.r.t.
any complex structure
compatible with the euclidean metric.}
\end{remark}

\medskip

{\scriptsize
\noindent {\sc
Claudio Arezzo}\\
Dipartimento di Matematica,
Universit\`a di Parma\\
Via M. D'Azeglio 85,
43100 Parma\\
{\tt claudio.arezzo@unipr.it}
}

\smallskip
{\scriptsize
\noindent {\sc Gian Pietro Pirola}\\
Dipartimento di Matematica,
Universit\`a di Pavia\\
Via Ferrata 1,
27100 Pavia\\
{\tt pirola@dimat.unipv.it}
}

\smallskip
{\scriptsize
\noindent {\sc Margherita Solci}\\
Dipartimento di Matematica,
Universit\`a di Pavia\\
Via Ferrata 1,
27100 Pavia\\
{\tt marghe@dimat.unipv.it}
}

\end{document}